\documentclass[12pt,thmsa]{article}
\usepackage{amsfonts}
\usepackage{amssymb}
\newtheorem{teo}{Theorem}[section]

\newtheorem{cor}[teo]{Corollary}

\newtheorem{obs2}[teo]{Remark}

\newtheorem{tea}{Theorem}[subsection]

\newtheorem{no2}[teo]{Note}

\newtheorem{no3}[tea]{Note}

\newcommand{\Gal}{{\rm Gal}}

\newcommand{\GL}{{\rm GL}}

\newcommand{\F}{{\mathbb{F}}}
\newcommand{\Q}{{\mathbb{Q}}}

\begin{document}
\title{{How to facet a gemstone: from potential modularity to the proof of Serre's modularity conjecture
}}
\author{Luis Dieulefait
\\
Dept. d'Algebra i Geometria, Universitat de Barcelona;\\
Gran Via de les Corts Catalanes 585;
08007 - Barcelona; Spain.\\
e-mail: ldieulefait@ub.edu\\
 }
\date{\empty}

\maketitle

\vskip -20mm

\begin{abstract}
In this survey paper we present recent results obtained by Khare, Wintenberger and the author that have led to a proof of Serre's conjecture, such as existence of compatible families, modular upper bounds for universal deformation rings and existence of minimal lifts, prime switching and modularity propagation, weight reduction (via existence of conjugates) and (iterated) killing ramification.
 The main tools used in the proof of these results are modularity lifting theorems \`{a} la Wiles and a result of potential modularity due to R. Taylor.

\end{abstract}

\section{Introduction}
We will present a series of results that emerged during the last five years and led to a proof of Serre's modularity conjecture. These results, obtained by Khare, Wintenberger and, independently, by the author, use as main tools the potential modularity result proved by R. Taylor in [Ta1] and [Ta2] together with modularity lifting theorems \`{a} la Wiles. The idea is that one should try to show modularity of a given representation by ``connecting" it somehow to a modular one. So it is a game where each move consists of ``switching" to a different characteristic where the situation becomes simpler. The two main theorems needed to make this possible are the ``existence of families" result in [Di1] and the ``existence of minimal lifts" in [Di2] and [K-W1]. Taylor's potential modularity is the crucial input allowing to prove both of these results.\\
 With these two theorems proved, combining them with modularity lifting theorems, one can start to make moves seeking to propagate modularity to more and more general situations (or, if you see it the other way round, to reduce the proof of modularity to a more and more simple situation): in fact it is like a fire which starts in very simple cases of representations with small level, weight and characteristic where modularity was known since the '70, and then with the two new results we propagate modularity from these base cases to all possible cases. The technical conditions needed to apply modularity lifting theorems make necessary some detours, were if not for these technicalities the proof of Serre's conjecture for any level and weight would follow automatically from a  base case by just a naive set of moves to reduce the level to $1$ (by moving to the primes in the level, one after the other) and then to switch to the case of characteristic $3$, where modularity of level $1$ representations was known since 1973. \\
 So, in this survey paper we will present the key ideas in the proof of the above two main results and we will explain in detail how they suffice to give a proof of the level $1$ weight $2$ case of Serre's conjecture. We will also show how to obtain a proof of the level $1$ case and arbitrary weight by a nice argument allowing to ``connect" a case of weight $k$ to a case of weight approximately $k/2$ via a new trick called ``existence of Galois conjugates" and some elementary considerations on Galois actions on roots of unity. Finally, we will explain how using the new and powerful modularity lifting results of Kisin, the proof of the odd conductor case can be reduced to the proof of the level $1$ case by performing ``iterated killing ramification". To overcome the technical conditions needed to apply this result of Kisin one has to perform some specific moves, following a combination of ideas of Khare-Wintenberger (to deal with the condition on the size of the residual image) and the author (to deal with local conditions).\\
 At the end of this section we also include, for the reader's convenience, a prelude to recall some concrete examples of modularity lifting theorems and how they work.\\
 
 Acknowledgments: This survey is based on two talks I gave last summer, at the Segundas Jornadas de Teor\'{i}a de N\'{u}meros (Madrid) and at the Summer School on Serre's Modularity Conjecture (Luminy), respectively. I would like to thank A. Quir\'{o}s, J. Cilleruelo, C. Corrales and E. Gonzalez, organizers of the SJTN, and J.-M. Fontaine, P. Colmez, M. Harris, R. Taylor and M. Kisin, organizers of the SSSMC, for the invitation.

\subsection{Prelude: How to apply Modularity Lifting Theorems?}

Modularity lifting theorems that initiated with the celebrated work of Wiles and Taylor-Wiles on Fermat's Last Theorem (cf. [W], [T-W]) will be one of the key tools that we will need.
By combining different modularity lifting results, we observed already in 2002 (cf. [D-M]) that in the case of crystalline $p$-adic Galois representations of small weights a fairly general modularity lifting theorem can be deduced:

\begin{teo} 
\label{teo:crista}
Let
$$         \rho_p : G_\Q  \rightarrow   \GL_2 (\bar{\Q}_p) $$
be a continuous, odd and irreducible Galois representation, with finite ramification set. Suppose that $\rho_p$ is crystalline, of
 Hodge-Tate weights $\{  0, k-1 \}$, $k \geq 2$. If $p> k$, $p \neq 2k-3$, and the residual representation $\bar{\rho}_p$ is either modular or reducible, then $\rho_p$ is modular.
\end{teo}

Starting from this result, we obtained the following (cf. [D-M]):

\begin{teo} \label{teo:calabi}
Let $X$ be a rigid Calabi-Yau threefold defined over $\Q$ having good reduction at $3$ and $7$. Then $X$ is modular, i.e., the compatible family of $2$-dimensional Galois representations $\{ \rho_\ell \}$ corresponding to the action of the group $G_\Q$ on the $\ell$-adic cohomology groups $H^3_{et} (X_{\bar{\Q}}, \Q_\ell)$ is modular. 
\end{teo}

Remark: Since it follows from the definition of $X$ that the Hodge-Tate numbers of this family are $\{0 , 3 \}$, the corresponding modular form $f$ has weight $4$.\\

Tools for the proof of ``Theorem \ref{teo:crista} $\Rightarrow$ Theorem \ref{teo:calabi}": Taking $p=7$ the $7$-adic representation is crystalline of weights $\{0, 3 \}$ and theorem \ref{teo:crista} can be applied. The modularity (if irreducible) of the residual representation follows from previous results of Manoharmayum on Serre's conjecture over $\F_7$.\\

Remark: Using modularity lifting theorems of Kisin for the case of crystalline representations of intermediate weights, and results of Berger-Li-Zhu for such representations, theorem \ref{teo:crista} extends to the case $p = k-1$ and therefore it also follows the stronger result that all rigid Calabi-Yau threefolds defined over $\Q$ with good reduction at $3$ are modular (as observed by the author in 2005 in a personal communication to K. Hulek).\\

Idea of the proof of Theorem \ref{teo:crista}:
The idea is to show that under the running hypothesis, either $\rho_p$ is ordinary or the residual image is ``large", namely, if $K$ is the quadratic field unramified outside $p$, the restriction $\bar{\rho}_p |_{G_K}$ is absolutely irreducible. In the ordinary case, modularity of $\rho_p$ follows from the assumption that the residual representation is either reducible or modular by results of Skinner-Wiles, and in the other case, modularity follows from residual modularity by results of Diamond-Flach-Guo and Taylor.\\

To guarantee that we are in one of these two cases, we first studied the dihedral case following the description given by Serre and Ribet in their work on the determination of images of families of Galois representations. This way we observe that the dihedral case induced from $K$ can only occur if $p$ and $k$ satisfy some specific relations. Finally, we apply results of Breuil and M\'{e}zard that imply that under the running hypothesis if the residual representation is ordinary, so is $\rho_p$. As a conclusion, the only problematic case is the dihedral induced from $K$ such that the action of inertia on the residual representation is given by fundamental characters of order $2$, and this case can only occur if $p = 2k-3$.\\

\section{Serre's conjecture in cases of small level and weight}

We start with a Galois representations with values in a finite field, namely:

$$ \bar{\rho} : G_\Q  \rightarrow  \GL_2 (\bar{\F}_p)$$
an irreducible, continuous, odd, Galois representation with (finite) ramification set $S \cup \{ p \}$, where $p$ is an odd prime.
Note that the representation is required to be odd: this means that the determinant of the image of complex conjugation is $-1$.\\
Let $\F$ be the ``field of coefficients" of $\bar{\rho}$, i.e., the finite extension of $\F_p$ generated by the traces of all matrices in its image. It is easy to see that we can assume without loss of generality that the representation is defined over $\F$:
$$ \bar{\rho}: G_\Q  \rightarrow  \GL_2 (\F) $$
Serre defines two invariants of such a Galois representation (cf. [Se]) known as the Serre's level $N(\bar{\rho})$ and the Serre's weight $k(\bar{\rho})$. The first one has to do with ramification at the primes in the set $S$: it is the prime-to-$p$ part of the Artin conductor of $\bar{\rho}$. The Serre's weight is defined in terms of the ramification at $p$: one has to look at the image of the inertia group at $p$, it can be described by fundamental characters of level $1$ or $2$, and the definition of $k(\bar{\rho})$ depends mainly on the exponents these characters are raised to (the Fontaine-Laffaille weights of the representation, cf. [Se]). For example, in the case of $k(\bar{\rho}) = 2$ the description of the ramification at $p$ corresponds to the case of Fontaine-Laffaille weights $0$ and $1$, thus we expect a global lift corresponding to a crystalline representation of weights $0$ and $1$ (a Barsotti-Tate representation). In fact, one of the consequences of Taylor's potential modularity result (a result that we will discuss later), combined with Ramakrishna's construction of ``geometric" (in the sense of Fontaine-Mazur) lifts, is that if $k(\bar{\rho}) = 2$ then there exists a lift corresponding to a $\GL_2$-type abelian variety with good reduction at $p$.\\
The case of $k(\bar{\rho}) = p+1$ corresponds to semistable ramification (in the sense of Fontaine) at $p$, and it appears when we consider the residual representation attached to an abelian variety of $\GL_2$-type with semistable reduction at $p$ (except in cases of loose of ramification when even if the variety has semistable reduction the residual weight turns out to be $2$, as in Wiles' proof of Fermat's Last Theorem if we assume that we are in Sophie Germain's second case, i.e., that the exponent $p$ divides $a, b$ or $c$ in  $a^p + b^p = c^p$).\\
Serre's conjecture (cf. [Se]) predicts that for this kind of residual representation $\rho$ there exists a modular form $f$ and a prime $P \mid p$ in the ring of integers of $\Q_f$ such that if $\rho_{f, P}$ is the Galois representation associated to $f$ by Shimura and Deligne then when reducing mod $P$ it holds:
$$ \bar{\rho}_{f, P} \cong \bar{\rho} $$
In other words, that $\bar{\rho}$ is always modular, in the sense that it admits a modular $p$-adic lift.\\
There is also a strong version of Serre's conjecture where the level and weight of one of the modular forms $f$ lifting $\bar{\rho}$ is specified. Yes, your guess is right. The invariants $N(\bar{\rho})$ and $k(\bar{\rho})$ that we introduced above should agree with the level and weight of some modular lift (cf. [Se]).\\
The first spectacular result in the direction of Serre's conjecture is the proof that the weak version (asking only for a modular lift) and the strong version are equivalent. The proof is due to several authors, the most important results are the weight control obtained by B. Edixhoven (cf. [Ed]) and the results of lowering the level of K. Ribet (cf. [Ri1]). It is this celebrated result of Ribet the one that implies that Fermat's Last Theorem follows from the Taniyama-Shimura conjecture.\\
Remark: In a few exceptional cases the conjecture has to be modified and one has to consider mod $p$ modular forms as defined by Katz. \\

In this section we will restrict to the case $k(\bar{\rho}) =2$, and we also assume that the representation is semistable at each of the primes in $S$. In this case the determinant is just $\chi$, the mod $p$ cyclotomic character, and the level $N(\bar{\rho})$ is square-free. More precisely, for every $q$ in $S$ (if any) we are assuming that the image of inertia at $q$ is unipotent:\\
$$ \bar{\rho}|_{I_q} \cong 
	    \left(
	             \begin{array}{cc}
	                    1  &  *  \\
	                    0  &  1   \\
	              \end{array}
	    \right) 
	     $$

The goal of this section is to describe the proof of a few cases of Serre's conjecture, cases where both the level and weight of $\bar{\rho}$ are small. This proof was obtained by the author and, independently, by Khare-Wintenberger, during 2004 (cf. [Di2], [K-W1]). For simplicity, we take the case of a representation $\bar{\rho}$ with Serre's level $1$ and weight $2$.\\
Let us briefly recall the main two tools that we will use in the proof:\\
(1) Modularity lifting theorems: theorem \ref{teo:crista} will be enough for this section. For more general cases of the proof of Serre's conjecture, more modularity lifting theorems are needed. All these results are of the form: Given a $p$-adic Galois representation $\rho$ such that the residual representation $\bar{\rho}$ is modular or reducible, then assuming some conditions on $\rho$ (mainly local conditions) and on $\bar{\rho}$ (mainly on its image) we can conclude that $\rho$ is modular. \\
Modularity lifting theorems have been proved by Taylor-Wiles, Diamond, Skinner-Wiles, Diamond-Flach-Guo, Fujiwara, Taylor, Savitt, and more recently the more general and powerful of them are due to M. Kisin.\\
(2) Potential Modularity (R. Taylor, cf. [Ta1], [Ta2]): Given a $p$-adic representation $\rho$
$$         \rho : G_\Q  \rightarrow   \GL_2 (\bar{\Q}_p) $$
which is odd, irreducible, continuous, with finite ramification set, crystalline of Hodge-Tate weights $\{0, k-1 \}$ with $p > k$ and $p \neq 2k-3$ (this last condition can be replaced by the condition: image of $\bar{\rho}$ non-solvable), then there exists a totally real Galois number field $F$ such that $p$ is unramified in $F/\Q$ (and has residual degree at most $2$) where the restriction of $\rho$ is modular, namely, $\rho|_{G_F}$ agrees with the $p$-adic Galois representation attached to a Hilbert cuspidal modular form over $F$, of parallel weight $k$.  \\
Remark: the result can also be applied to a residual representation with non-solvable image, for example one can first construct a suitable $p$-adic lift by results of Ramakrishna and then apply to it the result of Taylor. \\

Observe that in the level $1$ weight $2$ case of Serre's conjecture, modularity is equivalent to non-existence. In fact, since the weak and the strong version of Serre's conjecture are equivalent, in this case modularity gives the existence of a lift corresponding to a cuspidal modular form of level $1$ and weight $2$, but since there are no such modular forms, we deduce non-existence of the residual representation.\\

The starting point is the fact that Serre's conjecture is true in this case if $p=3$. The non-existence of irreducible representations of Serre's level $1$ in characteristic $3$ was proved by Serre in 1973, following a similar result of Tate for $p=2$. Their proofs use bounds for discriminants. \\

In order to prove this case in higher characteristic, the argument is the following: reduce the proof to the case of $p=3$ solved by Tate and Serre. This will be possible once we prove a result of ``switching  the residual characteristic". \\
In order to do so, we will have to apply the results (1) and (2), and we have to prove the two main results:\\
(i) Existence of strictly compatible families (cf. [Di1])\\
(ii) Existence of minimal lifts (cf. [Di2] and [K-W1])\\
In particular, the result (i) combined with modularity lifting theorems and the result of Tate-Serre gives a proof of the case of  conductor $1$ and ``weight $2$" (i.e., crystalline of weights $0$ and $1$) of the Fontaine-Mazur conjecture (cf. [Di1]).\\

With these two main results and modularity lifting theorems the method of ``switching the residual characteristic" (an idea that the author and J.-P. Wintenberger had independently) is the following: Let $\bar{\rho}$ be a representation of Serre's level $1$ and weight $2$ of characteristic $p > 3$. Using (ii), we construct a minimal $p$-adic lift $\rho$, minimality implies that its conductor  is $1$ and it is crystalline of weights $0$ and $1$. Using (i) we insert $\rho$ in a strictly compatible family $\{ \rho_\lambda \} $, thus all members of this family have conductor $1$ and are crystalline of the same weights: $0$ and $1$. Now we look at a $3$-adic member of this family, let us call it $\rho_3$ for simplicity. If we consider the residual representation $\bar{\rho}_3$ since it has level $1$ the result of Tate-Serre implies that it is reducible. Then we can apply a modularity lifting result of Skinner and Wiles, because the residual characteristic is $3$ and the weights of the lift are $0$ and $k-1 = 1$ (thus $3 > k = 2$) and then 
 this is contained in the more general modularity lifting theorem \ref{teo:crista}. The conclusion is that $\rho_3$ is modular, thus by compatibility $\rho$ is also modular, and $\bar{\rho}$ also by definition. This finishes the proof. As we already stressed, modularity here implies non-existence.\\
 
 So, what remains now is to explain how to prove the main results  (i) and (ii). We will see that in both cases the key tool to obtain them is the potential modularity result of Taylor (2). \\
 
\subsection{Existence of strictly compatible families (results of Dieulefait-Taylor)}

For simplicity we stick to the ``weight $2$" case, i.e., we consider a $p$-adic representation $\rho$ (two-dimensional, odd) which is crystalline of Hodge-Tate weights $0$ and $1$. We can also assume that it has non-solvable residual image. It can be shown (using potential modularity) that the field of coefficients of $\rho$ is a number field $M$.
Potential modularity combined with solvable base change of Arthur and Clozel implies that the restriction of $\rho$ to the Galois group of certain totally real number field $F$ is modular and moreover that the same is true if we restrict $\rho$ to the Galois group of any $E \subseteq F$ such that $\Gal(F/E)$ is solvable. In particular, it follows from this that the restrictions of $\rho$ to a collection of subfields $E_i$ of $F$ belong to compatible families, because Galois representations attached to Hilbert modular forms always live in compatible families. Using Brauer's theorem one can write $\rho$ in terms of its restrictions to these fields $E_i$, and since these restrictions belong to compatible families, this gives us a way of defining for any other prime $p'$ a virtual $p'$-adic representation $\rho'$ which is, by construction, ``compatible" with $\rho$: we just write the same formula that expresses $\rho$ in terms of its modular restrictions to the fields $E_i$ but we change from the $p$-adic to the $p'$-adic member in each of the modular compatible families. It remains to see that this construction gives us a true Galois representation, not just a virtual one. To see this, one can compute its inner product and using the fact that $\rho$ is a true Galois representation one concludes using the compatibility in each of the modular families used to build $\rho'$: the key point is that because of this compatibility and Cebotarev density theorem, the condition for being a true Galois representation is seen to be ``independent of $p$", thus if it holds for one prime (and we know it does), it holds for all. \\
To conclude, we also want to show that the compatible family containing $\rho$ that we have built satisfies some expected good local properties. A consequence of potential modularity observed by Taylor is that for two different representations in a potentially modular compatible family, corresponding to two primes $p_1$ and $p_2$, the ramification behavior at any third prime $p_3$ is the same for both of them. We can also show that for any prime $q$ such that $\rho$ is unramified at $q$, the $q$-adic representation in the family is crystalline. This follows from potential modularity and known properties of Galois representations attached to Hilbert modular forms if we assume that $q$ is unramified in $F$, and in the general case one uses again solvable base change to reduce to a case where $q$ (or at least a place of $F$ dividing it) is unramified. 
 
\subsection{Existence of minimal lifts}  

Let $\F$ be the  field of coefficients of $\bar{\rho}$. We will assume that the image of $\bar{\rho}$ is non-solvable. In fact in the weight $2$ semistable case (with $p > 3$) it is enough to assume that $\bar{\rho}$ is irreducible because then results of Ribet (cf. [Ri2]) imply that the image is as large as possible, namely the group of all matrices in $\GL_2(\F)$ with determinant in
 $\F_p^*$. \\
 
 {\bf Boeckle's lower bound}:
 We start by recalling results of Boeckle on minimal deformations (cf. [Bo1], [Bo2]). We consider minimal deformations of $\bar{\rho}$, this is a global deformation problem with local conditions, the conditions being that the ramification set of the deformation is the same as the one of $\bar{\rho}$, that the deformations are also semistable at ramified primes different from $p$ and Barsotti-Tate at $p$. In particular we are forcing the deformations to have fixed determinant, which is the cyclotomic character $\chi$. As usual in deformation theory, coefficient rings of deformations are complete Noetherian local rings with residue field 
 $\F$. Let us call $R$ the corresponding universal object, i.e., the universal minimal deformation ring of $\bar{\rho}$. It is well-known that $R$ has the following form:
 $$ R \cong \frac{  W(\F) [[X_1,....,X_n]]   }{(f_1,.....,f_m)}  $$
 where $W(\F)$ is the ring of Witt vectors, and $n$ and $m$ are unknown.\\
 A result of Boeckle gives a lower bound for these universal rings. His proof uses a local-to-global principle and the  determination of the corresponding versal local deformation rings. The lower bound can be stated as $n \geq m$ in the formula above, equivalently, it says that the absolute Krull dimension of $R$ is greater than $0$: $\dim R > 0$.\\
 From this, Boeckle himself obtained as a corollary that it would be enough to obtain a suitable upper bound for $R$ in order to conclude the existence of minimal $p$-adic lifts of $\bar{\rho}$. More precisely, he proved the following:\\
 \begin{cor} \label{teo:bo} (Boeckle): Assume that $R$ is finite (i.e., finitely generated as a $W(\F)$-module). For example, assume that we have a surjection:
 $$ R' \rightarrow R $$
 for some finite $W(\F)$-algebra $R'$. Then $R$ is finite flat complete intersection and, in particular, there are minimal $p$-adic lifts of $\bar{\rho}$.
  \end{cor}
  
  {\bf Potential modularity and the modular upper bound for $R$ (results of Dieulefait/Khare)}: To complete Boeckle's results we use potential modularity to obtain a suitable upper bound. We know that the restriction of $\bar{\rho}$ to the Galois group of some totally real number field $F$ (such that $p$ is unramified in it) is modular. We base change to $F$ and we compare $R$ with its analog for the restriction $\bar{\rho}|_{G_F}$, let us call it $R'$. Because of potential modularity and modularity lifting theorems it is known that $R'$ agrees with a Hecke algebra $\mathbb{T}$, and in particular that it is  finite flat complete intersection (this is an analog over $F$ of the results of Taylor-Wiles). The existence of minimal modular deformations in the context of Hilbert modular forms, which generalizes Ribet's celebrated lowering-the-level result, was proved by Jarvis, Rajaei and Fujiwara. So it remains to explain how we compare the rings $R$ and $R'$ and deduce the modular upper bound that completes Boeckle's results. This can be done in two slightly different ways, two approaches that were independently deviced during 2004 by Khare and the author (respectively). In Khare's approach, it is shown that the fact that $R'$ is finite implies that $R$ is also finite (see [K-W1] for details), which is enough due to corollary \ref{teo:bo} to conclude existence of minimal lifts. In our approach, we first ``manipulate" a little bit Taylor's potential modularity result to show that there is a certain freedom in choosing the field $F$ of modularity, and if we make an optimal (in a sense to be explained below) choice of $F$ we do have a surjection from $R'$ to $R$, thus by corollary \ref{teo:bo} existence of minimal lifts follows. An optimal choice consists on choosing $F$ such that:\\
  (i) it is linearly disjoint from the field fixed by the kernel of $\bar{\rho}$, and\\
  (ii) the degree of $F$ over $\Q$ is prime to $p$.\\
  In [Di2], it is proved that the field $F$ where potential modularity is achieved can be taken to satisfy these two conditions, in the weight $2$ semistable case. First, the large image result of Ribet is shown to be enough to deduce condition (i). For condition (ii), one just takes any Galois field $F$ satisfying (i) and consider a $p$-Sylow subgroup of the Galois group of $F$ over $\Q$. Using once again solvable base change we also have potential modularity over the subfield of $F$ which is fixed by this $p$-Sylow, which of course satisfies condition (ii) and still satisfies condition (i). \\
  Remark: for a residual representation of any level and weight, with non-solvable image, it was later observed (2005) by Taylor, Harris and Shepherd-Barron and, independently, by the author, that one can modify the proof of potential modularity (by slightly improving the Moret-Bailly argument by imposing more local conditions) to show that one can always choose the field $F$ to be linearly disjoint from any previously given number field, thus in particular $F$ can always be chosen to satisfy condition (i) (and by solvable base change also (ii)).\\
  
  Finally, taking an $F$ that satisfies these two conditions we observe that since by condition (i) the image of $\bar{\rho}|_{G_F}$ is equal to the image of $\bar{\rho}$, by condition (ii) and the fact that the kernel of reduction mod $p$ is a pro-$p$ group we deduce that for any minimal deformation of
   $\bar{\rho}|_{G_F}$ that extends to a minimal deformation of $\bar{\rho}$ the image is not increased in the descent. Thus when we restrict the minimal universal deformation of $\bar{\rho}$ to the Galois group of $F$ its coefficient ring is still $R$, thus giving the desired surjection from $R'$ to $R$ (by the universal property).\\
   We conclude that starting from $R' = \mathbb{T}$  we got a modular upper bound for $R$, thus providing exactly the condition under which Boeckle's corollary \ref{teo:bo} holds. In particular existence of minimal $p$-adic lifts is proved.
  
   \section{The level 1 case}
   
The result in the previous section can also be used to show some other cases of Serre's conjecture of small level and weight. Using essentially the same strategy, one gets cases of weight $2$ and level $2, 3, 5, 7, 11$ or $13$, and also the cases of level $1$ and weights $2,4,6,8,12,14$ (cf. [K-W1] and [Di2]). Using weight $2$ lifts, and also some non-minimal lifts with specific ramification conditions, Khare was able to extend the above results in 2005 and prove by induction on the weight (or, equivalently, in the characteristic) the level $1$ case of Serre's conjecture for arbitrary weight $k$ (cf. [Kh]). Later (2007) we obtained a different proof based on a new result of ``existence of Galois conjugates" (cf. [Di4]). In this section we will explain this new proof. \\
 The cases of $k \leq 14$, $k \neq 10$ will be the base for a proof by induction on the weight. One main feature of the results of the previous section (generalized to the case of arbitrary $k$) is that once the principle of ``switching the residual characteristic" is established one can deduce from it that in order to prove the level $1$ case for some fixed weight $k$, if we assume that it holds for one prime $p_0 \geq k-1$, then it holds for every prime $p \geq k-1$. This follows by exactly the same argument used in the previous section in the weight $2$ case to translate the proof of Serre's conjecture for $p_0 = 3$ to a proof for any prime. Observe that the ``switching" only covers primes $p \geq k-1$ but this is not a serious restriction, because for a smaller prime $p'$ one can assume that the Serre weight is at most $p' + 1$ (after twisting we can always reduce to this case) so by the induction hypothesis we can assume that in characteristic $p'$ the level $1$ case is already solved.\\
 This being said, if we are given a weight $k$ either equal to $10$ or larger than $14$ our task is to prove the level $1$ case of Serre's conjecture for this weight in some well-chosen characteristic $p_0 \geq k-1$. To use the induction hypothesis we need to establish some bridge to connect a suitable lift of the given representation to some representation known-to-be-modular, for example, a representation with level $1$ and smaller weight $k' < k$, and we need the bridge to be modular-preserving, i.e., we need that modularity lifting theorems \`{a} la Wiles can be applied to propagate modularity. \\
 For example, a naive try is to just apply the method of the previous section in a case of weight $k$ arbitrary. We start in some characteristic $p$, we take a minimal lift, which is a $p$-adic representation unramified outside $p$ crystalline of weights $0$ and $k-1$, we insert it in a compatible family and we switch to the prime $3$, the corresponding mod $3$ representation is reducible by the result of Serre, but the problem is that modularity lifting results do not apply in this situation where the characteristic is smaller than $k-1$, as we stressed in theorem \ref{teo:crista} (recall that we remarked in section 1 that this theorem extends to the case of $p=k-1$). \\
 Instead of minimal lifts as defined in the previous section, we will consider minimal weight $2$ lifts, following Khare (cf. [Kh]). 
 Let $\bar{\rho}$ be a representation of weight $k >2$ and values on a finite field of characteristic $p$. Then we consider a deformation which is minimal outside $p$ and at $p$ is potentially Barsotti-Tate, with determinant $\chi \omega^{k-2}$ where $\chi$ is the $p$-adic cyclotomic character and $\omega$ is the Teichmuller lift of the mod $p$ cyclotomic character. Moreover we can assume that the local inertial parameter at $p$ is exactly $(\omega^{k-2} \oplus 1, 0)$. The proof of existence of this type of lifts follows the same pattern that the proof of existence of minimal lifts given in section 1: the local-to-global principle of Boeckle and computations of  dimensions of local versal deformation rings give also in this case a lower bound for the corresponding universal deformation ring, and potential modularity can be used to obtain by base changing to some totally real number field $F$ a complementary upper bound. Observe that in the modular world the existence of lifts of this type follows from the theory  of congruences between modular forms. The advantage of this weight $2$ lifts is that now the weight is encoded in the ``nebentypus"
  $\omega^{k-2}$, and we will manipulate this lift trying to make the weight change in a controlled way. \\
  To gain a new degree of freedom, we need to introduce a new result: ``existence of Galois conjugates". Again, this is a well-known result if we assume that we have a representation which is modular, i.e., attached to some modular form $f$: in this case the conjugate representation will be the representation attached to $f^\sigma$ for $\sigma$ an element of the Galois group of (the Galois closure of) the field generated by the eigenvalues of $f$. Given a $p$-adic $2$-dimensional Galois representation $\rho$ (odd, irreducible), we see that the same arguments used in the previous section to prove existence of families also give the following: If potential modularity can be applied to $\rho$, then for any element $\sigma$ of the Galois group of the normal closure of the field $E$ of coefficients of $\rho$, there exists a conjugate Galois representation $\rho^\sigma$ whose traces and ramification are related to those of $\rho$ as in the modular case (i.e., the case where $\rho$ corresponds to $f$ and $\rho^\sigma$ to $f^\sigma$). As in the proof of existence of families, the conjugation is first performed on the intermediate fields $E_i$ where we do have modularity, then we can define the conjugate as a virtual representation as in the previous section, and check that it is a true representation because $\rho$ is so.\\
  Now that we have all the tools we need, let us explain how the induction goes: If we start with a weight $k > 14$ or $k = 10$ we will choose as a residual characteristic $p$ the first prime greater than $k$, except if $k= 32$ (this case will be discussed later). So  we call $\bar{\rho}$ a level $1$ weight $k$ representation in this characteristic $p$, and we consider a weight $2$ lift $\rho$ of it. The field of coefficients of $\rho$ contains of course the field generated by the values of the nebentypus $\omega^{k-2}$, so if we conjugate by a suitable element $\sigma$ we can change the nebentypus at will. If we consider a conjugate representation $\rho^\sigma$ it is clear that the modularity of $\rho$ is equivalent to the modularity of $\rho^\sigma$. \\
  The goal is to find a good $\sigma$, such that when we consider the residual representation $\bar{\rho^\sigma}$, also of level $1$, using the fact that the inertial parameter at $p$ of $\rho^\sigma$ is $(( \omega^{k-2})^\sigma \oplus 1, 0)$ when we compute the possible values of the Serre weight (as usual, after twisting so that they become at most $p+1$) of  $\bar{\rho^\sigma}$ they are smaller than $k$. Modularity lifting theorems (and in particular, the results of Skinner-Wiles in the reducible case) do apply in these potentially Barsotti-Tate situations (by a generalization of theorem \ref{teo:crista}   that covers this case) and we can conclude, if we find such a $\sigma$, that if $\bar{\rho^\sigma}$ is modular or reducible, then $\bar{\rho}$ is modular, thus finishing the proof of the level $1$ case by induction. \\
   {\bf Existence of a good conjugation $\sigma$:}  For small weights $k=10$ and $k > 14$ up to $k=36$ the existence of $\sigma$ as above can be checked by hand. We will list examples of this later, just observe that in the exceptional case $k=32$ the initial characteristic $p$ should be taken to be $p = 43$.\\
   For $k > 36$ the existence of $\sigma$ follows easily from the following estimates: If $p_n$, $p_{n+1}$ are consecutive primes and $p_{n+1} > 37 $, then:
   $$ \frac{p_{n+1} - 1}{ p_n  - 1}  < 1.15 < 1.2 \qquad \quad (2.1) $$
   
   which follows from known estimates on primes. \\
   So if $d= (k-2, p-1)$ where $p$ is the first prime larger than $k$, denoting $m = (p-1)/d$, $\omega^{k-2}$ and $(\omega^{k-2})^\sigma$ have values which are $m$-roots of unity and by (2.1) we know that $m>6$ if $k>36$. For smaller values of $k$ (either $k=10$ or $k>14$) we check by hand that $m>6$ except for $m=5$ in the case $k=10, p=11$. \\
   Since $m>6$ or $m=5$ we define $\sigma$ such that:
   $$ (\omega^{k-2})^{\sigma} = (\omega^d)^t $$
   for an exponent $t$ prime to $m$ and ``close to" $m/2$. More precisely, we take $t$ to be a number relatively prime to $m$ such that its distance to $m/2$ is at most $2$ (we take it bigger than $m/2$), and we know that it will not be equal to $m-1$ because $m=5$ or $m>6$. \\
   With this definition, we can compute the two possible values of the Serre's weight of $\bar{\rho^\sigma}$ (as usual, up to twist)
    and  we obtain the values $k_1 = dt+2$ and $k_2 = p+3- k_1$ which are both near $p/2$ for this $t$, thus since $k$ was larger than the prime previous to $p$ and using again estimates on primes it is easy to conclude that in particular we obtain $k_1 , k_2 < k$. We leave as an exercise to the reader to check in detail that this inequality holds for all $k > 36$ (or see [Di4] for details).\\
    Let us show the computations for a couple of small values of $k$ starting with the exceptional case $k=32$. Here we take $p=43$. We have $k-2 = 30$, $p-1 = 42$, $d=6$, $m = 7 > 6 $. We take $t= 4$. After conjugating we obtain a representation $\rho^\sigma$ with nebentypus $\omega^{dt} = \omega^{24}$. Observe that we started with a $\rho$ with nebentypus $\omega^{k-2} = \omega^{30}$. Computing the possible Serre's weight of this conjugate representation we obtain the values $ k_1 = d t+ 2 = 26$ and $k_2 = p+3 - k_1 = 20$, both smaller than $k=32$, so the induction works even in this exceptional case. Observe that we were forced to take $p=43$ because for $p = 37$ or $p=41$ we end up with $m=6$ and $m= 4$ (respectively) and there is no way to choose a $t$ that makes the induction work, the only conjugated nebentypus that we can get for such an $m$ is by taking $\sigma$ equal to complex conjugation, and this may not change the Serre's weight.\\
    Here are four other small values of $k$ where we check at hand that the induction works:\\
    $$ k=10, p= 11, d= 2, m= 5, t= 3, k_i= 6, 8 $$
    $$ k=20, p= 23, d= 2, m= 11, t= 6, k_i= 12, 14 $$
    $$ k=24, p= 29, d= 2, m= 14, t= 9, k_i= 12, 20 $$
    $$ k=34, p= 37, d= 4, m= 9, t= 5, k_i= 18, 22 $$

   \section{The odd level case}
   As we mentioned in the introduction, the ideal proof would be the following: if $q_1, ...q_r$ are the primes in the level $N$, by ``switching the prime" to $q_1$, then to $q_2$, and so on, we finally reduce the proof to the level $1$ case, already solved. The problem is that modularity lifting theorems have technical conditions. So we have to do better. Luckily, by the beginning of 2006, Kisin presented a new and very powerful modularity lifting result (cf. [Ki]), a result that in particular worked for liftings of non-zero Hodge-Tate weight $k-1$ with $k$ arbitrary large (recall the usual restriction to small $k$, as in theorem \ref{teo:crista}). With this result we have shown how to reduce the proof of Serre's conjecture for any case of odd level $N$ to the level $3$ case (cf. [Di3]), a case that can be solved as the level $1$ case, using as base cases some results of Schoof of modularity of abelian varieties with small odd conductor. The reduction to the level $3$ case is similar to the above naive process of iterated killing ramification, except that first we have to make some moves to reduce to a situation where the technical conditions needed to apply the main theorem of [Ki] are satisfied. The technical conditions are of two types:\\
   (a) residual images must be non-solvable\\
   (b) a local condition on the residual representation, and the lift required to be potentially semistable or crystalline over an abelian extension\\
   (a) is dealt with by applying a trick of Khare-Wintenberger (cf. [K-W2]): they introduce a ``good-dihedral prime", one can reduce to a case where there is such a prime in the level, and the ramification at this prime forces the residual images to be non-solvable.\\
   (b) this is dealt with in [Di3] by, via weight $2$ lifts, changing the set of primes $q_1,.....,q_r$ in the level to another set
   $q'_1, ..... q'_r$ where ramification is ``potentially crystalline over a cyclic extension of small degree" and the primes are ``far from each other", in such a way that the local technical conditions needed to apply Kisin's result are satisfied, and then the iterated killing ramification as described before (just ``switch" to each of the primes in the level except the good-dihedral prime, one after the other, in increasing order) can be performed.\\
   Remark: The moves that have to be done to go from the original primes $q_i$ to the $q'_i$ are of two types: the first set of moves change via weight $2$ lifts to pseudo Sophie-Germain primes $q''_i$ where (by using as ``pivots" the two odd prime divisors
    of $q''_i-1$) one is reduced to a situation of semistable ramification. Then, in the second set of moves, one changes to primes 
    $q'_i$ which are such that $q'_i \equiv 1 \; \pmod{q''_i - 1}$.\\
    
    In the above process, because of the technical condition (a) we are assuming that there is a good-dihedral prime $q$ in the level as in [K-W2] that guarantees that the residual images are large. At the end, one has to remove also this extra prime from the level, but it is convenient before doing so to first reduce the weight as much as possible, and it turns out that $k=4$ is the minimal weight that we can obtain here (cf. [Di3]). By switching to characteristic $3 = k-1$ and taking a semistable weight $2$ lift, i.e., corresponding to an abelian variety with semistable reduction at $3$, we are now ready to kill part of the ramification at the good-dihedral prime as in [K-W2] (by switching to the odd prime $t$ dividing the order of the ramification at $q$), and we end up with a case of weight $2$ and level $3q$, and by a standard argument this is equivalent to a weight $k= q+1$, level $3$ case (just switch to characteristic $q$ and use the fact that for an abelian variety with semistable reduction at $q$ the mod $q$ Serre's weight is $q+1$ or $2$). We have thus reduced the proof to the level $3$ case. \\
    Applying again weight reduction it is enough to deal with the cases $k=4,6$ and level $3$ ($k=2$ been already proved), which by an easy trick are reduced to the case of $k=2$ and level $21$. Let us explain this step in more detail: we use the Sophie Germain pair of primes $7$ and $3$. First we move to characteristic $7$ where the weight is transformed into nebentypus $\omega^{k-2}$ as in the previous section, then we switch to the prime $3$ because $3$ is the order of $\omega^{k-2}$, thus in the mod $3$ representation the only ramification at $7$ that can appear is semistable ramification.
    This means that we can now take a weight $2$ lift of this residual representation having semistable ramification at $3$ and $7$, which is known to correspond (by Taylor's potential modularity result) to an abelian variety with this ramification. Since R. Schoof proved that for such a small odd conductor every semistable abelian variety is modular, the proof of the level $3$ case is concluded.\\
    
    Remark: In May 2006 we devised this strategy to prove the odd conductor case of Serre's conjecture using [Ki] as main tool, and at the same time Khare-Wintenberger where completing a different proof of the odd case of Serre's conjecture which is independent of [Ki]: their proof uses previous modularity lifting results and it consists on a double induction on the level and the weight, so you need first to reduce the weight $k$ to $2$, then by killing ramification, i.e., switching to a prime in the level, you make the level smaller but the weight increases again, then again do weight reduction until arriving to $k=2$, and so on. The problem is that in every weight reduction they do need to go down to weight $k=2$ for their method to work, and in order to do so they need to work at several steps in characteristic $2$. On the one hand, their strategy was deviced  before May 2006 (this is the content of [K-W2]), but on the other hand, by May 2006 the modularity lifting results in characteristic $2$ that they needed where not available. Nevertheless shortly after they proved such results in a second part to their previous work (cf. [K-W3]), thus making unconditional the proof given in [K-W2]. Notice finally that together with the results of [Ki] and new results of Schoof, my proof incorporates the key trick of introducing a good-dihedral prime in the level, which is taken from [K-W2]. \\
    Let us mention to conclude that Kisin has recently obtained more general modularity lifting results in characteristic $2$ that allow to cover also the remaining cases of Serre's conjecture, i.e., the cases of even level.\\

\section{Bibliography}

[Bo1] Boeckle, G., {\it 
A local-to-global principle for deformations of Galois representations},
J. Reine Angew. Math. {\bf 509} (1999) 199-236 
\newline
[Bo2] Boeckle, G, {\it On the isomorphism $R_\emptyset \rightarrow T_\emptyset$}, appendix to {\it On isomorphisms between deformation and Hecke rings} by C. Khare.,
Inventiones Math. {\bf 154} (2003)  217-222 
\newline
[Di1] Dieulefait, L., {\it Existence of compatible families and new cases of the Fontaine-Mazur conjecture},
 J. Reine Angew. Math. {\bf 577} (2004) 147-152 
 \newline
 [Di2] Dieulefait, L., {\it The level $1$ weight $2$ case of Serre's conjecture}, preprint (2004), to be published by December 2007 
\newline
[Di3] Dieulefait, L., {\it Remarks on Serre's modularity conjecture}, preprint (2006)
\newline
[Di4] Dieulefait, L., {\it The level $1$ case of Serre's conjecture revisited}, preprint (2007)
\newline
[D-M] Dieulefait, L., Manoharmayum, J., {\it Modularity of rigid Calabi-Yau threefolds over $\Q$}, in
``Calabi-Yau Varieties and Mirror Symmetry", Fields Institute Communications Series, AMS {\bf 38} (2003) 159-166
\newline
[Ed] Edixhoven, B., {\it The weight in Serre's conjectures on modular forms}, Invent. Math. {\bf 109} (1992) 563-594
\newline
[Kh] Khare, C., {\it Serre's modularity conjecture: the level $1$ case}, Duke Math. J. {\bf 134} (2006)  557-589
\newline
[K-W1] Khare, C., Wintenberger, J.-P., {\it On Serre's conjecture for $2$-dimensional mod $p$ representations of the Galois group of $\Q$}, preprint (2004)
\newline
[K-W2] Khare, C., Wintenberger, J.-P., {\it Serre's modularity conjecture (1)}, preprint (2006)
\newline
[K-W3] Khare, C., Wintenberger, J.-P., {\it Serre's modularity conjecture (2)}, preprint (2006)
\newline
[Ki] Kisin, M., {\it The Fontaine-Mazur conjecture for $\GL_2$}, preprint (2006)
\newline
[Ri1] Ribet, K., {\it On modular representations of $\Gal(\bar{\Q} / \Q)$ arising from modular forms}, Invent. Math. {\bf 100} (1990) 431-476
\newline 
[Ri2] Ribet, K., {\it Images of semistable Galois
representations}, Pacific J.  Math. {\bf 181} (1997)
\newline
[Se] Serre, J-P., {\it Sur les repr{\'e}sentations modulaires de degr{\'e}
$2$ de $\Gal(\bar{\mathbb{Q}} / \mathbb{Q})$}, Duke Math. J. {\bf 54}
(1987) 179-230
\newline
%[SeIII] Serre, J-P., {Oeuvres}, Vol. III (1986) Springer-Verlag
%\newline
[Ta1] Taylor, R., {\it Remarks on a conjecture of Fontaine and Mazur},
J. Inst. Math. Jussieu {\bf 1} (2002)
\newline
[Ta2] Taylor, R., {\it On the meromorphic continuation of degree two
 L-functions}, Documenta Mathematica, Extra Volume: John Coates' Sixtieth Birthday (2006) 729-779 
\newline
[T-W] Taylor, R., Wiles, A., {\it Ring theoretic properties of certain Hecke algebras},
Ann. of Math. {\bf 141} (1995) 553-572
\newline
[W] Wiles, A., {\it Modular Elliptic Curves and Fermat's Last Theorem}, Ann.
of Math. {\bf 141} (1995) 443-551

\end{document}